\documentclass[preprint,12pt]{elsarticle}

\usepackage{amssymb}

\journal{Applied Numerical Mathematics }

\usepackage{amsmath}

\usepackage{graphicx}
\usepackage{theorem}
\usepackage{amssymb}
\usepackage{color}

\theorembodyfont{\upshape}

\newtheorem{theorem}{Theorem}[section]
\newtheorem{lemma}[theorem]{Lemma}

\begin{document}

\begin{frontmatter}
\title{Two-grid economical algorithms for parabolic integro-differential equations with nonlinear memory
}


\author{Wansheng Wang\fnref{wsw}}
\ead{w.s.wang@163.com}

\author{Qingguo Hong\corref{cor1}\fnref{qgh}}
\ead{huq11@psu.edu}

\cortext[cor1]{Corresponding author}

 \address[wsw]{Department of Mathematics, Shanghai Normal University, Shanghai, 200234,
China}
\address[qgh]{Department of Mathematics,
        Pennsylvania State University, State College, PA 16802,U.S.A.}

\begin{abstract}
In this paper, several two-grid finite element algorithms for solving parabolic integro-differential
equations (PIDEs) with nonlinear memory are presented. Analysis of these
algorithms is given assuming a fully implicit time discretization. It is shown that these algorithms are as stable as the standard fully discrete finite element algorithm, and can achieve the same accuracy as the standard algorithm if the coarse grid size $H$ and the fine grid size $h$ satisfy $H=O(h^{\frac{r-1}{r}})$.
Especially for PIDEs with nonlinear memory defined by a lower order
nonlinear operator, our two-grid algorithm can save significant
storage and computing time. Numerical experiments are given to confirm the theoretical results.
\end{abstract}


%

\end{frontmatter}

\section{Introduction}
\label{intro} The main purpose of this paper is to present some discretization
techniques based on two finite element subspaces for solving
parabolic integro-differential equations (PIDEs) with nonlinear
memory:
\begin{eqnarray}\label{eqNV1}
  &&u_t+Au+\int^t_{0}K(t-s)Bu(s)ds =f(x,t),~~(x,t)\in \Omega\times (0,T],\\
\label{intcon1} &&u(x,t)=0,~~~~(x,t)\in \partial \Omega \times (0,T],\\
\label{intcon2}&&u(x,0)=u_0(x),~~~~x\in \Omega,\end{eqnarray} where
$\Omega\subset {\mathbb R}^d (d\ge 1)$ is a bounded and polyhedral
domain with a piecewise smooth boundary $\partial \Omega$, $K(t)$ is
a smooth or nonsmooth memory kernel, and $f$ is a known function, $u$ is
the solution we need to solve which is scalar function .
$A$ is a symmetric positive definite second-order elliptic operator
with smooth coefficients in $x$ and $t$, and $B$ is a nonlinear
operator of at most second order; that is,
\begin{eqnarray}\label{eq1.4}
Bu=-\nabla \cdot \left(\alpha(x,u)\nabla u+\beta (x,u)\right)+\gamma (x,u)\cdot \nabla u+g(x,u).
 \end{eqnarray}
 For brevity, we will drop the dependence of variable $x$ in $\alpha(x,u)$, $\beta(x,u)$, $\gamma(x,u)$, and $g(x,u)$ in the following exposition.
We assume that the functions $\alpha(u)\in {\mathbb R}^{d\times d}$ is a tensor function, $\beta(u))\in {\mathbb R}^d,~\gamma(u)\in {\mathbb R}^d$ are vector functions, and $g(u)\in {\mathbb R}^1$ is scalar function, respectively. And all the functions $\alpha(u)$, $\beta(u)$, $\gamma(u)$ and $g(u)$ are smooth and bounded together with the Gateaux derivative.
 For the functions $\beta(u)$ and $g(u)$, we also assume that $\beta(0)=0$ and $g(0)=0$.

Equations of the above type, or linear versions thereof, can arise
from many physical processes in which it is necessary to take into
account the effects of memory due to the deficiency of the usual
diffusion equations \cite{Gurtin68,Miller78,Raynal79}. For
approximating the solution $u$ of PIDEs, both finite difference and
finite element methods have been investigated extensively in the
past for both the linear and nonlinear problem (see, for example,
\cite{Cannon88,Cannon90,Lin90,Lopez-Marcos90,Pani96,Chen98,Zhang09}).
Recently, several new numerical methods such as mixed finite element
method, finite volume element method, and discontinuous Galerkin
method for space discretization or time discretization have been
proposed to solve PIDEs (see, for example,
\cite{Ewing02,Sinha06,Pani08,Sinha09,Biswas10,Mustapha11}).

The two-grid method based on two finite element spaces, one on a
coarse grid and one on a fine grid, was first developed by Xu
\cite{Xu921,Xu922,Xu94,Xu96} for nonsymmetric linear and nonlinear
elliptic problems. Since then, the two-grid method for elliptic
problems has been investigated further, e.g., Axelsson and Layton
\cite{Axelsson96}, Xu and Zhou \cite{Xu01}, Li and Huang
\cite{Li10}, and Bi and Ginting \cite{Bi07,Bi11}. In these works,
theoretical study and numerical experiments show that the combined
use of the numerical method such as finite element method and finite
difference method, and the two-grid technique is computationally
more efficient than the original method. Due to this better
practical performance, the two-grid method has been widely applied
to the study of eigenvalue problems \cite{Xu012,Xu02,Hu11}, steady
Navier Stokes equations \cite{Layton98,Girault01,He05,Dai08}, the
time-dependent Navier Stokes problem
\cite{Girault012,Abboud08,Abboud09,Shang10,Tachim08}, the nonlinear
parabolic problem
\cite{Dawson94,Marion95,Dawson98,Chen10,Wu99,Chen03,Chen09,Qin05},
and nonlinear hyperbolic equations \cite{Chen103}. Recently, Jin,
Shu, and Xu \cite{Jin06} used this technique to solve decoupling
systems of partial differential equations; Mu and Xu \cite{Mu07} and
Cai, Mu, and Xu \cite{Cai09} employed it for the mixed Stokes-Darcy
model. In \cite{Wang}, we proposed the two-grid algorithms
based on the backward Euler scheme and finite element approximation
for semi-linear PIDEs, and studied the long-time stability and error
estimates of the two-grid algorithms.

In this paper, we present some two-grid algorithms for PIDEs with
nonlinear memory and perform theoretical analysis that demonstrates
our methods' ability to match the accuracy of the classic finite
element method by (1) solving a nonlinear problem
on a coarse space $S_H$ and (2) solving a symmetric positive definite linear
problem on the fine space $S_h$. Thus, solving PIDEs with nonlinear memory
is not much more difficult than solving one linear problem, as $dim
~S_H\ll dim ~S_h$ and the work involved in solving the nonlinear
problem on the coarse grid is relatively limited.

It is worth adding that when $\alpha\equiv0$, our algorithm
significantly reduces computational memory and storage requirements.
A practical difficulty of numerical methods for PIDEs is that all
previous values must be stored, as they all enter subsequent
equations. In order to reduce memory requirements, some economical
schemes have been proposed (for example, see,
\cite{Sloan86,Huang94}). However, these schemes either require more
regularities on the solution $u$ \cite{Sloan86}, or they cannot be
applied to nonlinear problem \cite{Huang94}.

The remainder of this article is organized as follows: In Section 2,
we present some conventions and notations that will be used
throughout the article. In Section 3, the stability and error
estimate of the classic fully discrete finite element method are
discussed. The two-grid algorithms for PIDEs with nonlinear memory
are presented and the stability and error estimates of these
algorithms are discussed in Section 4. In Section 5, we offer some
concluding remarks.

Throughout this paper, we use the letters $C$ and $c$ (with and
without subscripts) to denote a generic positive constant that stand
for different values depending on the context in different
equations. When it is not important to keep track of these
constants, we conceal the letter $C$ or $c$ in the notation
$\lesssim$ or $\gtrsim$, such that $x\lesssim y$ means $x\le Cy$ and
$x\gtrsim y$ means $x\ge cy$.

\section{Preliminaries}
\setcounter{equation}{0}
For any non-negative integer $r$ and number $p\ge 1$, let $\mathcal
W^{r,p}(\Omega)$ be the standard Sobolev space with a norm
$\|\cdot\|_{r,p}$ given by $\|v\|^p_{r,p}=\sum_{|\kappa|\le
r}\|D^\kappa v\|^p_{\mathcal L^p(\Omega)}$ (with the usual
modification if $p=\infty$). This Sobolev space is also equipped
with the seminorm $|v|^p_{r,p}=\sum_{|\kappa|= r}\|D^\kappa
v\|^p_{\mathcal L^p(\Omega)}$. For $p=2$, we denote $\mathcal
H^r=\mathcal W^{r,2}(\Omega)$ and take $\mathcal H^1_0$ as the
subspace of $\mathcal H^1$ consisting of functions with a vanishing
trace on $\partial \Omega$. For simplicity, we also use notations
$\|\cdot\|_r$, $\|\cdot\|$ and $\|\cdot\|_\infty$, and $|\cdot|_r$
such that $\|\cdot\|_r=\|\cdot\|_{r,2}$, $\|\cdot\|=\|\cdot\|_{0,2}$
and $\|\cdot\|_\infty=\|\cdot\|_{0,\infty}$, and
$|\cdot|_r=|\cdot|_{r,2}$.

Let $\{S_h\}_{0<h\le 1}$ be a family of finite-dimensional subspaces
of $\mathcal H^1_0$, with the following approximation properties:
\begin{eqnarray}\label{eq2.1} \inf\limits_{\chi\in S_h}\left\{\|v-\chi\|+h\|v-\chi\|_1\right\}
\lesssim h^r\|v\|_r,\qquad v\in \mathcal H^r\cap \mathcal H^1_0,\quad r\ge 1+\frac{d}{2}.
\end{eqnarray}
We also assume that $\{S_h\}_{0<h\le 1}$ satisfies the inverse
hypothesis: there exists a constant $C>0$ independent of $h$ such
that
\begin{eqnarray}\label{eq2.2}
\|\nabla \chi\|_\infty\le Ch^{-d/2}\|\nabla \chi\|, \qquad \chi\in
S_h.
\end{eqnarray}

The weak formulation of the problem (\ref{eqNV1}), (\ref{intcon2})
is: Find $u\in \mathcal H^1_0(\Omega)$ such that
\begin{eqnarray}\label{eq2.3}
&&(u_t,v)+A(u,v)+\int^t_{0}K(t-s)B(u(s),v)ds =(f,v),~~v\in H^1_0,\\
&&u(0)=u_0,\label{eq2.4} \end{eqnarray} where $A(\cdot,\cdot)$ is
the bilinear form associated with the operator $A$ on $\mathcal
H^1_0\times \mathcal H^1_0$ and $B(\cdot,\cdot)$ is defined by
$$B(u,v)=(\alpha(u)\nabla u+\beta(u),\nabla v)+(\gamma(u)\cdot \nabla u+g(u),v), \qquad u,v\in \mathcal W^{1,\infty}\cap \mathcal H^1_0.$$
$(\cdot,\cdot)$ denotes the inner product in $\mathcal L^2(\Omega)$.
We always assume that $A$ is coercive and continuous with coercivity constant $\nu_0$ and continuity constant $\nu_1$. That is, we have
\begin{eqnarray}\label{eq2.5}
A(v,v)&\ge& \nu_0\|v\|^2_1\qquad \forall v\in\mathcal
H^1_0,\\
\label{eqa2.5}
|A(u,v)|&\le& \nu_1\|u\|_1\|v\|_1\qquad \forall u, v\in\mathcal
H^1_0.\end{eqnarray} In view of the assumptions on the functions
$\alpha(u),~\beta(u),~\gamma(u)$, and $g(u)$, it is easily verified
that there exists a positive constant $\mu_0$ such that
\begin{eqnarray}\label{eq2.6}
|B(u,v)|\le \mu_0\|u\|_1\|v\|_1.
\end{eqnarray}

For the time discretization of (1.1)-(1.3) we will consider the
backward Euler scheme. To analyze the discretization on a time
interval $(0,T]$, let $N$ be a positive integer, $\Delta t=T/N$, and
let $t_n=n\Delta t$. As the truncation error of the backward Euler
scheme is $O(\Delta t)$, we introduce a quadrature formula with a
truncation error $O(\Delta t)$,
\begin{eqnarray}\label{eq2.7} \Delta t\sum\limits_{i=1}^n\omega_{ni}g(t_i)=\int^{t_n}_0K(t_n-s)g(s)ds+O(\Delta t).
\end{eqnarray}
Given our emphasis on two-grid discretization in space, we will not
discuss how to obtain the numbers $\omega_{ni}$, but only assume
that there exists a positive constant $K_1$ such that
$|\omega_{ni}|\le K_1$ for any $1\le n\le N,~1\le i\le n$ and that
$\omega_{nn}\not=0$. Therefore, the problem considered in this paper
must be discretized by a fully implicit scheme. Thus, the backward
Euler fully discrete finite element approximation of problem (1.1),
(1.3) is defined as a sequence $\{U^n\}_{n=0}^N$, such that
\begin{eqnarray}\label{eq2.8}
&&\left(\bar \partial U^{n},v\right)+A(U^{n},v)+\Delta t\sum\limits_{i=1}^n\omega_{ni} B(U^i,v) =(f^n,v),~v\in S_h,~ n\ge 1,\\
&&U^0=u^h_0,\label{eq2.9} \end{eqnarray} where $\bar \partial U^{n}=\frac{U^{n}-U^{n-1}}{\Delta t}$, $u^h_0$ is an
appropriate approximation of $u_0$ in $S_h$, $f^n=f(t_n)$. We know
that (\ref{eq2.8}) will result in a truncation error $O(\Delta t)$
in time. But for nonlinear problems considered in this paper
($\omega_{nn}\not=0$), the solution of a nonlinear algebraic system
is required at each time step. To decrease the amount of
computational work, we propose using a two-grid technique to solve
the PIDEs with nonlinear memory. With this technique, at each time
step, solving a nonlinear problem on the fine space $S_h$ is reduced
by solving a nonlinear problem on the coarse space $S_H$ and solving a
linear SPD problem on the fine space $S_h$.

For functions that vanish on the boundary, we recall Poincare's inequality: there exists a constant $\mathcal P$ such that
$$\forall v\in \mathcal H^1_0, \qquad \|v\|\le \mathcal P|v|_1.$$

We make extensive use of the $\epsilon-$type inequality $2ab\le
\epsilon a^2+b^2/\epsilon,~\epsilon>0$, and of the inequality
$a^2+b^2\le (|a|+|b|)^2$. The results of this paper are based on the
identity
\begin{eqnarray}\label{eq2.10}
2(a^{n+1},a^{n+1}-a^n)=|a^{n+1}|^2-|a^{n}|^2+|a^{n+1}-a^n|^2,
\end{eqnarray}
and the following Gronwall lemma proved in \cite{Emmrich05}.

\begin{lemma}[Discrete Gronwall lemma \cite{Emmrich05}] Let $0\le \lambda<1$, and $a_n,~b_n,~c_n,~\lambda_n\ge 0$ with $\{c_n\}$ being monotonically increasing. Then \begin{eqnarray}
a_n+b_n\le \sum\limits_{j=\varpi}^{n-1}\lambda_ja_j+\lambda a_n+ c_n,\qquad n=\varpi,\varpi+1,\cdots
\end{eqnarray}
implies for $n=\varpi,\varpi+1,\cdots$
\begin{eqnarray*}
a_n+b_n\le \frac{c_n}{1-\lambda}\prod\limits_{j=\varpi}^{n-1}\left(1+\frac{\lambda_j}{1-\lambda}\right)\le \frac{c_n}{1-\lambda}\exp\left(\frac{1}{1-\lambda}\sum\limits_{j=\varpi}^{n-1}\lambda_j\right).
\end{eqnarray*}
\end{lemma}

\section{Error estimate for the classic fully discrete finite element method}\setcounter{equation}{0}
In this section, we discuss the stability and error estimate of the
standard fully discrete finite element method (\ref{eq2.8}),
(\ref{eq2.9}). First, we prove the stability of the solution of
(\ref{eq2.8}) and (\ref{eq2.9}).

\begin{theorem} \label{th3.1}
Let $U^n$ be the solution obtained by (\ref{eq2.8}) and
(\ref{eq2.9}). Then for all \begin{eqnarray} \label{eqr3.1}\Delta
t\le\min\left\{\frac{1}{2},\frac{7\nu^2_0}{8\mu_0^2K^2_1T}\right\},\end{eqnarray}
we have
\begin{eqnarray}\label{eq3.1}
&&\|U^n\|+\left(\sum\limits_{i=1}^n\|U^i-U^{i-1}\|^2\right)^{1/2}
+\frac{\sqrt{\nu_0}}{2} \left(\sum\limits_{i=1}^n\Delta t\|U^i\|^2_1\right)^{1/2}\nonumber\\
&\le& E^{1/2}_n\left(\|U^0\|^2+\Delta t \sum\limits_{i=1}^n
\|f^i\|^2 \right)^{1/2},
\end{eqnarray}
where $E_n=6\max\{e^{2t_n},e^{(2\mu_0K_1t_n/\nu_0)^2}\}$.
\end{theorem}

{\bf Proof.}
By taking $v=2\Delta t U^{n}$ in (\ref{eq2.8}) and using
(\ref{eq2.10}), we obtain
\begin{eqnarray}\label{eq3.2}
&&\|U^n\|^2-\|U^{n-1}\|^2+\|U^n-U^{n-1}\|^2+2\nu_0\Delta t \|U^n\|^2_1+2(\Delta t)^2 \sum\limits_{i=1}^n\omega_{ni}B(U^i,U^n)\nonumber\\
&\le& 2\Delta t \|f^n\|\|U^n\|.
\end{eqnarray}
Using (\ref{eq2.6}), we have
\begin{eqnarray}\label{eq3.3}
&&\|U^n\|^2-\|U^{n-1}\|^2+\|U^n-U^{n-1}\|^2+2\nu_0\Delta t \|U^n\|^2_1\nonumber\\
&\le& \mu_0(\Delta t)^2
\sum\limits_{i=1}^{n}|\omega_{ni}|(\frac{1}{\epsilon}\|U^i\|^2_1+\epsilon\|U^n\|^2_1)
+\Delta t \left(\|U^n\|^2+\|f^n\|^2\right).
\end{eqnarray}
Choose $\epsilon=\nu_0/(\mu_0K_1t_n)$ to obtain
\begin{eqnarray}\label{eq3.4}
&&\|U^n\|^2+\|U^n-U^{n-1}\|^2+\nu_0\Delta t \|U^n\|^2_1\nonumber\\
&\le& \|U^{n-1}\|^2+\frac{\mu_0^2K_1^2t_n}{\nu_0}(\Delta t)^2
\sum\limits_{i=1}^{n}\|U^i\|^2_1+\Delta t \|f^n\|^2+\Delta t
\|U^n\|^2.
\end{eqnarray}
By summation, we have
\begin{eqnarray}\label{eq3.5}
&&\|U^n\|^2+\sum\limits_{i=1}^n\|U^i-U^{i-1}\|^2+\nu_0 \Delta t \sum\limits_{i=1}^n\|U^i\|^2_1\nonumber\\
&\le& \|U^0\|^2+\Delta t \sum\limits_{i=1}^{n}\|U^i\|^2+(\Delta t)^2
\sum\limits_{i=1}^{n}\frac{\mu_0^2K_1^2t_i}{\nu_0}\sum\limits_{j=1}^i
 \|U^j\|^2_1+\Delta t \sum\limits_{i=1}^n\|f^i\|^2, \end{eqnarray}
which implies that
\begin{eqnarray}\label{eq3.6}
&&(1-\Delta t)\|U^n\|^2+\sum\limits_{i=1}^n\|U^i-U^{i-1}\|^2+\left(\nu_0-\frac{\mu_0^2K_1^2t_n}{\nu_0}\Delta t\right)
\Delta t \sum\limits_{i=1}^n\|U^i\|^2_1\nonumber\\
&\le& \|U^0\|^2+\Delta t \sum\limits_{i=1}^{n-1}
\max\left\{1,\frac{4\mu_0^2K_1^2t_i}{\nu_0^2}\right\}\left(\|U^i\|^2+\frac{\nu_0}{4}\Delta
t\sum\limits_{j=1}^i \|U^j\|^2_1\right) +\Delta t
\sum\limits_{i=1}^n\|f^i\|^2.\nonumber\\
\end{eqnarray} Since
condition (\ref{eqr3.1}) implies that $1-\Delta t\ge \frac{1}{2}$
and $\nu_0-\frac{\mu_0^2K_1^2t_n}{\nu_0}\Delta t\ge \frac{\nu_0}{8}$,
with the aid of discrete Gronwall lemma 2.1, we obtain
\begin{eqnarray}\label{eq3.7}
&&\|U^n\|^2+2\sum\limits_{i=1}^n\|U^i-U^{i-1}\|^2+\frac{\Delta t \nu_0}{4} \sum\limits_{i=1}^n\|U^i\|^2_1\nonumber\\
&\le& E_n\left(\|U^0\|^2+\Delta t \sum\limits_{i=1}^n \|f^i\|^2
\right),
\end{eqnarray}
which implies (\ref{eq3.1}). Thus the proof is completed.\quad

{\it Remark.} From (\ref{eqr3.1}), we find that for a given integral interval $(0,T]$ the stepsize $\Delta t$ is determined by the ratio of $\nu_0$ to $\mu_0$ and increases as the value of coercivity constant $\nu_0$ increases.

Now let us take $v=2\Delta t \bar \partial U^{n}$ in (\ref{eq2.8}) to obtain
\begin{eqnarray}\label{eqa3.2}
&&2\Delta t\|\bar \partial
U^n\|^2+A(U^n,U^n)-A(U^{n-1},U^{n-1})+A(U^n-U^{n-1},U^n-U^{n-1})\nonumber\\
&&+2(\Delta t)^2 \sum\limits_{i=1}^n\omega_{ni}B(U^i,\bar \partial
U^n)\nonumber\\ &=& 2\Delta t (f^n,\bar \partial
U^n).
\end{eqnarray}
Since
\begin{eqnarray}
2\Delta t (f^n,\bar \partial
U^n)&\le&\frac{1}{2}\Delta t\|f^n\|^2+2\Delta t\|\bar \partial
U^n\|^2
\end{eqnarray}
and
\begin{eqnarray}
2(\Delta t)^2 \sum\limits_{i=1}^n|\omega_{ni}B(U^i,\bar \partial U^n)|&\le&
2\Delta t\mu_0\sum\limits_{i=1}^n|\omega_{ni}| \|U^i\|_1\|U^n-U^{n-1}\|_1\nonumber\\
&\le& \frac{t_n\mu_0^2K_1^2}{\nu_0}\Delta t\sum\limits_{i=1}^n \|U^i\|^2_1+\nu_0\|U^n-U^{n-1}\|^2_1,
\end{eqnarray}
(\ref{eqa3.2}) becomes
\begin{eqnarray}\label{eqr83.2}
   \nu_0\|U^n\|^2_1 &\le& \frac{1}{2}\Delta t\|f^n\|^2+\frac{t_n\mu_0^2K_1^2}{\nu_0}\Delta t\sum\limits_{i=1}^n \|U^i\|^2_1+\nu_1\|U^{n-1}\|^2_1.
\end{eqnarray}
Then we have the following result.

\begin{theorem} \label{thr83.1}
Let $U^n$ be the solution obtained by (\ref{eq2.8}) and
(\ref{eq2.9}). Then for all \begin{eqnarray} \label{eqra3.1}\Delta
t\le \frac{\nu^2_0}{2\mu_0^2K^2_1T},\end{eqnarray}
we have
\begin{eqnarray}\label{eqr83.1}
\|U^n\|_1\le& C\left(\|U^0\|^2_1+\Delta t \sum\limits_{i=1}^n
\|f^i\|^2 \right)^{1/2}.
\end{eqnarray}
\end{theorem}
{\bf Proof.} It follows from (\ref{eqra3.1}) that $\nu_0-\frac{\mu_0^2K_1^2t_n}{\nu_0}\Delta t\ge \frac{\nu_0}{2}$. Then an application of discrete Gronwall lemma 2.1 to (\ref{eqr83.2}) leads to
$$\|U^n\|^2_1 \le C\left(\|U^0\|^2_1+\Delta t \sum\limits_{i=1}^n
\|f^i\|^2\right), $$
which implies (\ref{eqr83.1}). This completes the proof.\quad

To estimate the error of the fully discrete approximation
(\ref{eq2.8}), we define, for
$w,~u,~v\in \mathcal W^{1,\infty}\cap \mathcal H^1_0(\Omega)$,
$$B_1(w;u,v)=(\alpha(w)\nabla u,\nabla v)+(\gamma(w)\cdot \nabla u,v).$$
Due to the assumptions on $\alpha(u)$ and $\gamma(u)$, there exist a constant $\sigma$ such that
\begin{eqnarray}\label{eq3.8a}
|B_1(w;u,v)|\le \sigma\|u\|_1\|v\|_1.
\end{eqnarray}
As usual, we write the error $e^n=u(t_n)-U^n$ as
$$e^n=u(t_n)-U^n=(u(t_n)-V_hu(t_n))+(V_hu(t_n)-U^n)=\rho^n+\theta^n,$$
where $V_hu$ is the Ritz-Volterra projection of the solution $u$ and defined by \cite{Cannon90}
\begin{eqnarray}\label{eq3.8}
  A(u-V_hu,v)+\int^t_{0}K(t-s)B_1(u(s);u(s)-V_hu(s),v)ds =0,~~v\in S_h.
  \end{eqnarray}For $\rho(t)=u(t)-V_hu(t)$, following the line of Cannon and Lin \cite{Cannon90}, we show that there exists $C_0>0$, independent of $h$ and $t$,
  such that (see, also, \cite{Chen98,Lin91})
\begin{eqnarray}
  \|\rho(t)\|+h\|\rho(t)\|_1&\le& C_0h^r|||u(t)|||_r,\qquad t\ge 0,\label{eq3.9}\\
  \|\rho_t(t)\|&\le& C_0h^r\left(|||u(t)|||_r+|||u_t(t)|||_r\right),\label{eq3.10}\\
  \|\rho(t)\|_\infty&\le& C_0 h^r|\ln h||||u(t)|||_{r,\infty},\label{eq3.11}
  \end{eqnarray}
  where
  $$|||u(t)|||_r=\|u(t)\|_r+\int_0^t\|u(\tau)\|_rd\tau,\quad |||u(t)|||_{r,\infty}=\|u(t)\|_{r,\infty}+\int_0^t\|u(\tau)\|_{r,\infty}d\tau,$$
and there exists a positive
constant $C=C(u)$, independent of $h$, such that
\begin{eqnarray}\label{eq3.20a}
\|\nabla V_hu\|_\infty+\|\nabla (V_hu)_t\|_\infty\le C.
\end{eqnarray}

Now we need to estimate the error $\theta^n=V_hu(t_n)-U^n$.

\begin{theorem} \label{th3.2} Let $u$ and $U^n$ be the solutions of (\ref{eq2.3})-(\ref{eq2.4}) and (\ref{eq2.8})-(\ref{eq2.9}), respectively.
 If \begin{eqnarray}\label{eq3.12a}
 16\sigma^2K_1^2T< \nu_0^2,\end{eqnarray}
 then, for sufficiently small $\Delta t$, we have
\begin{eqnarray}\label{eq3.12}
\|\theta^n\|+\|\theta^n-\theta^{n-1}\|+\sqrt{\nu_0\Delta t}\|\theta^n\|_1\lesssim h^r+\Delta t.
\end{eqnarray}
\end{theorem}

{\bf Proof.} Firstly, it follows from (\ref{eq2.3}) and (\ref{eq2.8}) that
\begin{eqnarray*}
  &&\left(u_t-\bar\partial U^n,v\right)+A(u-U^n,v)\nonumber\\
  &&+\int^{t_n}_{0}K(t-s)B(u(s),v)ds-\Delta t\sum\limits_{i=1}^n\omega_{ni}B(U^i,v)=0,\qquad v\in S_h.
\end{eqnarray*}
Then we find that $\theta^n$ satisfies\begin{eqnarray}\label{eq3.13}
  &&\left(\bar\partial \theta^n,v\right)+A(\theta^n,v)+A(\rho^n,v)+\int^{t_n}_{0}K(t-s)B_1(u(s);u(s)-V_hu(s),v)ds\nonumber\\
  &&+\int^{t_n}_{0}K(t-s)B_1(u(s);V_hu(s),v)ds\nonumber\\
  &&+\int^{t_n}_{0}K(t-s)[(\beta(u(s)),\nabla v)+(g(u(s)),v)]ds-\Delta t\sum\limits_{i=1}^n\omega_{ni}B(U^i,v)\nonumber\\
  &=&-\left(\frac{\rho^n-\rho^{n-1}}{\Delta t},v\right)-\left(u_t-\frac{u(t_n)-u(t_{n-1})}{\Delta t},v\right),\qquad v\in S_h.
\end{eqnarray}
Using (\ref{eq3.8}), we have
\begin{eqnarray}\label{eqa1}
  &&\left(\bar\partial \theta^n,v\right)+A(\theta^n,v)+\int^{t_n}_{0}K(t-s)B_1(u(s);V_hu(s),v)ds\nonumber\\
  &&-\Delta t\sum\limits_{i=1}^n\omega_{ni}B_1(u(t_i);V_hu(t_i),v)+\Delta t\sum\limits_{i=1}^n\omega_{ni}B_1(u(t_i);\theta^i,v)\nonumber\\
  &&+\Delta t\sum\limits_{i=1}^n\omega_{ni}\left[\left(\left(\alpha(u(t_i))-\alpha(U^i)\right)\nabla U^i,\nabla v\right)+\left(\left(\gamma(u(t_i))-\gamma(U^i)\right)\cdot\nabla U^i,v\right)\right]\nonumber\\
  &&+\int^{t_n}_{0}K(t-s)[(\beta(u(s)),\nabla v)+(g(u(s)),v)]ds\nonumber\\
  &&-\Delta t\sum\limits_{i=1}^n\omega_{ni}\left[\left(\beta(u(t_i)),\nabla v\right)+\left(g(u(t_i)),v\right)\right]\nonumber\\
  &&+\Delta t\sum\limits_{i=1}^n\omega_{ni}\left[\left(\beta(u(t_i))-\beta(U^i),\nabla v\right)+\left(g(u(t_i))-g(U^i),v\right)\right]\nonumber\\
  &=&-\left(\frac{\rho^n-\rho^{n-1}}{\Delta t},v\right)-\left(u_t-\frac{u(t_n)-u(t_{n-1})}{\Delta t},v\right),\qquad v\in S_h.
\end{eqnarray}
Now, in view of (\ref{eq2.7}), we have
\begin{eqnarray}\label{eqa2}
\left|\int^{t_n}_{0}K(t-s)B_1(u(s);V_hu(s),v)ds-\Delta t\sum\limits_{i=1}^n\omega_{ni}B_1(u(t_i);V_hu(t_i),v)\right|\lesssim \Delta t\|v\|_1.\qquad
\end{eqnarray}
and
\begin{eqnarray}
&&\left|\int^{t_n}_{0}K(t-s)[(\beta(u(s)),\nabla v)+(g(u(s)),v)]ds\right.\nonumber\\
  &&\left.-\Delta t\sum\limits_{i=1}^n\omega_{ni}\left[\left(\beta(u(t_i)),\nabla v\right)+\left(g(u(t_i)),v\right)\right]\right|\lesssim \Delta t\|v\|_1.\qquad
\end{eqnarray}
Due to (\ref{eq3.8a}), the fifth term on the left-hand side in (\ref{eqa1}) can be bounded as
\begin{eqnarray}
\left|\Delta t\sum\limits_{i=1}^n\omega_{ni}B_1(u(t_i);\theta^i,v)\right|\le \Delta tK_1\sum\limits_{i=1}^n\sigma\|\theta^i\|_1\|v\|_1.
\end{eqnarray}
By virtue of the assumptions on $\alpha(u)$, $\beta(u)$, $\gamma(u)$ and $g(u)$, we know $\alpha$, $\beta$, $\gamma$ and $g$ satisfy Lipschitz conditions with Lipschitz constant $C_L$, and thus the sixth and ninth terms on the left-hand side in (\ref{eqa1}) are estimated as follows:
\begin{eqnarray}
&&\left|\Delta t\sum\limits_{i=1}^n\omega_{ni}\left[\left(\left(\alpha(u(t_i))-\alpha(U^i)\right)\nabla U^i,\nabla v\right)+\left(\left(\gamma(u(t_i))-\gamma(U^i)\right)\cdot\nabla U^i,v\right)\right]\right|\nonumber\\
&\le&\left|\Delta t\sum\limits_{i=1}^n\omega_{ni}\left[\left(\left(\alpha(u(t_i))-\alpha(U^i)\right)\nabla \theta^i,\nabla v\right)+\left(\left(\gamma(u(t_i))-\gamma(U^i)\right)\cdot\nabla \theta^i,v\right)\right]\right|\nonumber\\
&&\left|\Delta t\sum\limits_{i=1}^n\omega_{ni}\left[\left(\left(\alpha(u(t_i))-\alpha(U^i)\right)\nabla V_hu(t_i),\nabla v\right)\right.\right.\nonumber\\
&&\left.\left.+\left(\left(\gamma(u(t_i))-\gamma(U^i)\right)\cdot\nabla V_hu(t_i),v\right)\right]\right|\nonumber\\
&\le&\Delta tK_1\sum\limits_{i=1}^n\sigma\|\theta^i\|_1\|v\|_1+C_L\Delta tK_1\sum\limits_{i=1}^n\|u(t_i)-U^i\|\|\nabla V_hu(t_i)\|_\infty\|v\|_1\nonumber\\
&\le&\Delta tK_1\sum\limits_{i=1}^n\sigma\|\theta^i\|_1\|v\|_1+CC_L\Delta tK_1\sum\limits_{i=1}^n(\|\rho^i\|+\|\theta^i\|)\|v\|_1,
\end{eqnarray}
where the estimate (\ref{eq3.20a}) has been used, and
\begin{eqnarray}
&&\left|\Delta t\sum\limits_{i=1}^n\omega_{ni}\left[\left(\beta(u(t_i))-\beta(U^i),\nabla v\right)+\left(g(u(t_i))-g(U^i),v\right)\right]\right|\nonumber\\
&\le&C_L\Delta tK_1\sum\limits_{i=1}^n\|u(t_i)-U^i\|\|v\|_1\le C_L\Delta tK_1\sum\limits_{i=1}^n(\|\rho^i\|+\|\theta^i\|)\|v\|_1.
\end{eqnarray}
The first term on the right-hand side in (\ref{eqa1}) can be bounded as
\begin{eqnarray}\label{eq3.14}
    \left|\left(\frac{\rho^n-\rho^{n-1}}{\Delta t},v\right)\right|&\le& \frac{1}{\Delta t}\|\rho^n-\rho^{n-1}\|\|v\|\le\frac{1}{\Delta t}\int_{t_{n-1}}^{t_n}\|\rho_t(s)\|ds\|v\|\nonumber\\
     &\le& \frac{C_0 h^r}{\Delta t}\int_{t_{n-1}}^{t_n}\left(\||u(s)\||_r+\||u_t(s)\||_r\right)ds\|v\|;
  \end{eqnarray}
and the last term can be bounded as
\begin{eqnarray}\label{eqa8}\left|\left(u_t-\frac{u(t_n)-u(t_{n-1})}{\Delta t},v\right)\right|\le \int_{t_{j-1}}^{t_j}\|u_{tt}\|ds\|v\|.\end{eqnarray}

Taking $v=2\Delta t\theta^n$ and substituting all the above estimates (\ref{eqa2})-(\ref{eqa8}) into (\ref{eqa1}), we obtain
\begin{eqnarray}
&&\|\theta^n\|^2-\|\theta^{n-1}\|^2+\|\theta^n-\theta^{n-1}\|^2+2\nu_0\Delta t\|\theta^n\|_1^2\nonumber\\
&\le&C\Delta t^2\|\theta^n\|_1+C(\Delta t^2+\Delta th^r)\|\theta^n\|+4\sigma\Delta t^2K_1\sum\limits_{i=1}^n\|\theta^i\|_1\|\theta^n\|_1\nonumber\\
&&+4CC_L\Delta t^2K_1\sum\limits_{i=1}^n(\|\rho^i\|+\|\theta^i\|)\|\theta^n\|_1\nonumber\\
&\le&C\Delta t^2+C\Delta t^2\|\theta^n\|^2_1+C(\Delta t+h^r)^2+C\Delta t^2\|\theta^n\|^2+2\sigma\Delta t^2K_1\sum\limits_{i=1}^n\frac{1}{\epsilon_1}\|\theta^i\|^2_1\nonumber\\
&&+2\epsilon_1\sigma\Delta tK_1t_n\|\theta^n\|^2_1+2CC_L\Delta tK_1\sum\limits_{i=1}^n\|\rho^i\|^2+2CC_L\Delta t^2K_1t_n\|\theta^n\|_1^2\nonumber\\
&&+\frac{2}{\epsilon_2}CC_L\Delta t^2K_1\sum\limits_{i=1}^n\|\theta^i\|^2+2\epsilon_2CC_L\Delta tK_1t_n\|\theta^n\|^2_1,
\end{eqnarray}
where we have used the inequality
\begin{eqnarray*}
C\Delta t^2\sum\limits_{i=1}^na_ib_n&=&C\sum\limits_{i=1}^n\Delta t^{1/2}a_i\Delta t^{3/2}b_n
\nonumber\\
&\le& \frac{C}{2}\Delta t\sum\limits_{i=1}^na^2_i+\frac{C}{2}\Delta t^3\sum\limits_{i=1}^nb^2_n\le \frac{C}{2}\Delta t\sum\limits_{i=1}^na^2_i+\frac{C}{2}\Delta t^2t_n b^2_n.\end{eqnarray*}

Using the estimate (\ref{eq3.9}) for $\rho^i$, and taking $\epsilon_1=\frac{\nu_0}{4\sigma K_1t_n}$ and $\epsilon_2=\frac{\nu_0}{4CC_L K_1t_n}$, we have
\begin{eqnarray}
&&\|\theta^n\|^2-\|\theta^{n-1}\|^2+\|\theta^n-\theta^{n-1}\|^2+\nu_0\Delta t\|\theta^n\|_1^2\nonumber\\
&\le&C(\Delta t+h^r)^2+\left(C+\frac{8}{\nu_0}C^2C^2_LK_1^2t_n\right)\Delta t^2\|\theta^n\|^2\nonumber\\
&&+\left(C+\frac{8}{\nu_0}\sigma^2K^2_1t_n+2CC_LK_1t_n\right)\Delta t^2\|\theta^n\|^2_1+\frac{8}{\nu_0}\sigma^2\Delta t^2K^2_1t_n\sum\limits_{i=1}^{n-1}\|\theta^i\|^2_1\nonumber\\
&&+\frac{8}{\nu_0}C^2C^2_L\Delta t^2K^2_1t_n\sum\limits_{i=1}^{n-1}\|\theta^i\|^2,
\end{eqnarray}
Noting the condition (\ref{eq3.12a}) and taking sufficiently small $\Delta t$ such that $$\left(C+\frac{8}{\nu_0}\sigma^2K^2_1t_n+2CC_LK_1t_n\right)\Delta t\le \frac{\nu_0}{2}~~ {\hbox{and}}~~ \left(C+\frac{8}{\nu_0}C^2C^2_LK_1^2t_n\right)\Delta t\le \frac{1}{2},$$ we obtain
\begin{eqnarray}
&&\|\theta^n\|^2+\|\theta^n-\theta^{n-1}\|^2+\nu_0\Delta t\|\theta^n\|_1^2\nonumber\\
&\le&\|\theta^{n-1}\|^2+\frac{1}{2}\Delta t^2\sum\limits_{i=1}^{n-1}\|\theta^i\|^2+\frac{1}{2}\Delta t\|\theta^n\|^2+C(\Delta t+h^r)^2\nonumber\\
&&+\frac{\nu_0}{2}\Delta t\|\theta^n\|^2_1+\frac{\nu_0}{2}\Delta t^2\sum\limits_{i=1}^{n-1}\|\theta^i\|^2_1.
\end{eqnarray}

Applying discrete Gronwall lemma 2.1 to the above inequality yields
\begin{eqnarray}\label{eq3.15}
\|\theta^n\|^2+\|\theta^n-\theta^{n-1}\|^2+\nu_0\Delta t\|\theta^n\|^2_1\lesssim (h^r+\Delta t)^2.
\end{eqnarray}
which implies (\ref{eq3.12}). This proves the theorem.
  \qquad

Note that the condition (\ref{eq3.12a}), which implies that the equation (1.1) is diffusion-dominant, is appropriate, since the system may be blowup if the integral term is dominant. Under the condition (\ref{eq3.12a}), we can not study the long time behaviour of the numerical solution. Of course, if we assume that there exist positive constants $\alpha_0,~\alpha_1> 0$ such that
 \begin{eqnarray}\label{eq1.5}
  \alpha_0|\xi|^2\le \xi^{T}\alpha (u)\xi \le \alpha_1|\xi|^2, \quad \forall u\in \mathbb R, \quad \xi\in  \mathbb R^d,
   \end{eqnarray}
then following the approach of \cite{Wang}, we can study the long time behavior of the exact solution and the numerical solution.

We now give the $\mathcal H^1$ estimate of the error $\theta^n$.
\begin{theorem} \label{thr83.2} Let $u$ and $U^n$ be the solutions of (\ref{eq2.3})-(\ref{eq2.4}) and (\ref{eq2.8})-(\ref{eq2.9}), respectively. Then, for all $\Delta t$ satisfying
\begin{eqnarray}\label{eqt}
\Delta t< \frac{\nu_0^2}{16\sigma^2K_1^2T},
\end{eqnarray} we have
\begin{eqnarray}\label{eqr83.12}
\|\theta^n\|_1\lesssim h^{r}+\Delta t.
\end{eqnarray}
\end{theorem}

{\bf Proof.} Taking $v=2\Delta t\bar\partial \theta^n$ in (\ref{eqa1}), and estimating every terms in a way similar to Theorem 3.3, we get
\begin{eqnarray}
&&2\Delta t\|\bar\partial \theta^n\|^2+A(\theta^n,\theta^n)-A(\theta^{n-1},\theta^{n-1})+A(\theta^n-\theta^{n-1},\theta^n-\theta^{n-1})\nonumber\\
&\le&C\Delta t^2\|\bar\partial\theta^n\|_1+C(\Delta t^2+\Delta th^r)\|\bar\partial\theta^n\|+4\sigma\Delta t^2K_1\sum\limits_{i=1}^n\|\theta^i\|_1\|\bar\partial\theta^n\|_1\nonumber\\
&&+4CC_L\Delta t^2K_1\sum\limits_{i=1}^n(\|\rho^i\|+\|\theta^i\|)\|\bar\partial\theta^n\|_1\nonumber\\
&\le&\frac{C^2\Delta t^2}{\nu_0}+\frac{\nu_0}{4}\|\theta^n-\theta^{n-1}\|^2_1+\frac{C^2\Delta t}{8}(\Delta t+h^r)^2+2\Delta t\|\bar\partial\theta^n\|^2\nonumber\\
&&+\frac{16}{\nu_0}\sigma^2\Delta tK^2_1t_n\sum\limits_{i=1}^n\|\theta^i\|^2_1+\frac{\nu_0}{4}\|\theta^n-\theta^{n-1}\|^2_1+\frac{16}{\nu_0}C^2C^2_L\Delta tK^2_1t_n\sum\limits_{i=1}^n\|\rho^i\|^2\nonumber\\
&&+\frac{\nu_0}{4}\|\theta^n-\theta^{n-1}\|^2_1+\frac{16}{\nu_0}C^2C^2_L\Delta tK^2_1t_n\sum\limits_{i=1}^n\|\theta^i\|^2+\frac{\nu_0}{4}\|\theta^n-\theta^{n-1}\|^2_1.
\end{eqnarray}
Using (\ref{eq2.3}), (\ref{eq2.4}), (\ref{eq3.9}) and (\ref{eq3.15}) yields
\begin{eqnarray}
\nu_0\|\theta^n\|_1^2&\le&\nu_1\|\theta^{n-1}\|^2_1+C(\Delta t+h^r)^2+\frac{16}{\nu_0}\sigma^2\Delta tK^2_1t_n\sum\limits_{i=1}^n\|\theta^i\|^2_1.
\end{eqnarray}
Then when $\Delta t$ satisfies (\ref{eqt}), an application of discrete Gronwall lemma 2.1 to the above inequality leads to (\ref{eqr83.12}). This completes the proof
\quad

We observe that if (\ref{eq3.12a}) holds, then for any $\Delta t<1$, the conclusion (\ref{eqr83.12}) is valid.

In the next theorem, we will establish the error estimate for
the solution computed by the standard fully discrete finite element
method (\ref{eq2.8})-(\ref{eq2.9}).

\begin{theorem}[Error estimate for classic FEM] \label{th3.3}
Let $u$ be the solution of (\ref{eq2.3})-(\ref{eq2.4}) and $U^n$ be
the solution of (\ref{eq2.8})-(\ref{eq2.9}). Then, for sufficiently small $\Delta t$, we have, for all $n\ge 1$,
\begin{eqnarray}\label{eq3.16}
\|U^n-u(t_n)\|\lesssim  h^{r}+\Delta t,\qquad
\|U^n-u(t_n)\|_1\lesssim  h^{r-1}+\Delta t.
\end{eqnarray}

\end{theorem}

{\bf Proof.}The first inequality is a direct result of {\sc Theorem \ref{th3.2}} and (\ref{eq3.9}).
From {\sc Theorem \ref{thr83.2}} and (\ref{eq3.9}), we can prove the second inequality in (\ref{eq3.16}).\quad

\section{Two-grid algorithms for PIDEs with nonlinear memory}\setcounter{equation}{0}
In this section, we present three two-grid algorithms
of the backward Euler finite element method for PIDEs with nonlinear
memory. The basic mechanism in these algorithms is the construction
of two regular triangulations of $\Omega$: a coarse triangulation
$\mathcal T_H$ with mesh size $H$ and a fine one $\mathcal T_h$ with
mesh size $h$ ($h\ll H$). For practical purposes, $\mathcal T_h$ is
a refinement of $\mathcal T_H$. The corresponding finite element
spaces are $S_H$ and $S_h$, which will be called coarse and fine
space, respectively. To state the algorithms, we define, for
$w,~u,~v\in \mathcal W^{1,\infty}\cap \mathcal H^1_0(\Omega)$,
$$\tilde B(w;u,v)=(\alpha(w)\nabla u+\beta(w),\nabla v)+(\gamma(w)\cdot \nabla u+g(w),v).$$
Due to the assumptions on $\alpha(u),~\beta(u),~\gamma(u)$, and
$g(u)$, there exist two constants $\mu_1$ and $\mu_2$ such that
\begin{eqnarray}\label{eq4.1}
|\tilde B(w;u,v)|\le \mu_1\|u\|_1\|v\|_1+\mu_2\|w\|\|v\|_1.
\end{eqnarray}
Let us now present our first two-grid algorithm.

{\sc Algorithm 4.1.}

Step one (nonlinear problem on coarse grid $\mathcal T_H$): Given
$U^{n-1}_H$, find $U^{n}_H\in S_H$ such that
\begin{eqnarray}
  &&\frac{1}{\Delta t}\left(U^{n}_H-U^{n-1}_H,v\right)+A(U^{n}_H,v)+\Delta t\sum\limits_{i=1}^n\omega_{ni} \tilde B(U^i_H;U^i_H,v) =(f^n,v),\nonumber\\
  &&~~~~~~~~~~~~v\in S_H,~ n\ge 1,\label{eq4.2}\\
&&U^0_H=u^0_H,\label{eq4.3}\end{eqnarray}

Step two (linear problem on fine grid $\mathcal T_h$): Given
$U^{n}_H$, find $U^{n}_h \in S_h$ such that
\begin{eqnarray}
  &&\frac{1}{\Delta t}\left(U^{n}_h-U^{n-1}_h,v\right)+A(U^{n}_h,v)+\Delta t\sum\limits_{i=1}^n\omega_{ni} \tilde B(U^i_H;U^i_h,v)=(f^n,v),\nonumber\\
   &&~~~~~~~~v\in S_h,\qquad n\ge 1,\label{eq4.4}\\
&&U^0_h=u^0_h.\label{eq4.5}\end{eqnarray}

Firstly, we observe that for the solution of (\ref{eq4.4}) and
(\ref{eq4.5}), our stability result is similar to the solution of
(\ref{eq2.8}) and (\ref{eq2.9}).

\begin{theorem}[Stability of two-grid FEM  Algorithm 4.1] Let $U^n_h$ be the solution obtained by {\sc Algorithm 4.1}. Then
when $\Delta t$ satisfies (\ref{eqr3.1}) and \begin{eqnarray}\label{eqr4.6}\Delta
t\le\min\left\{\frac{1}{2},\frac{3\nu^2_0}{2\mu_1^2K^2_1T}\right\},\end{eqnarray}
we have
\begin{eqnarray}\label{eq4.6}
&&\sup\limits_{1\le i\le n}\|U^i_h\|+\left(\sum\limits_{i=1}^n\|U^i_h-U^{i-1}_h\|^2\right)^{1/2}+\frac{\sqrt{\nu_0}}{2}
\left(\sum\limits_{i=1}^n\Delta t\|U^i_h\|^2_1\right)^{1/2}\nonumber\\
&\le& C\left(\|U^0_h\|^2+\|U^0_H\|^2+\Delta t \sum\limits_{i=1}^n \|f^i\|^2 \right)^{1/2}.
\end{eqnarray}
\end{theorem}

{\bf Proof.} Similar to (\ref{eq3.3}), using (\ref{eq4.1}), we have
\begin{eqnarray}\label{eq4.7}
&&\|U^n_h\|^2-\|U^{n-1}_h\|^2+\|U^n_h-U^{n-1}_h\|^2+2\Delta t \nu_0\|U^n_h\|^2_1\nonumber\\
&\le& (\Delta t)^2
\sum\limits_{i=1}^{n}|\omega_{ni}|\left(\frac{\mu_1}{4\epsilon_1}\|U^i_h\|^2_1+\mu_1\epsilon_1\|U^n_h\|^2_1
+\frac{\mu_2}{4\epsilon_2}\|U^i_H\|^2+\mu_2\epsilon_2\|U^n_h\|^2_1\right)\nonumber\\
&&+\Delta t \left(\|U^n_h\|^2+\|f^n\|^2\right).
\end{eqnarray}
After choosing $\epsilon_1=\frac{\nu_0}{2\mu_1K_1t_n}$ and
$\epsilon_2=\frac{\nu_0}{2\mu_2K_1t_n}$, (\ref{eq4.7}) becomes
\begin{eqnarray}\label{eq4.8}
&&\|U^n_h\|^2+\|U^n_h-U^{n-1}_h\|^2+\Delta t \nu_0\|U^n_h\|^2_1\nonumber\\
&\le& \|U^{n-1}_h\|^2+(\Delta t)^2
\sum\limits_{i=1}^{n}\left(\frac{\mu^2_1K^2_1t_n}{2\nu_0}\|U^i_h\|^2_1+\frac{\mu^2_2K^2_1t_n}{2\nu_0}\|U^i_H\|^2\right)+\Delta
t \|f^n\|^2+\Delta t \|U^n_h\|^2.\nonumber\\
\end{eqnarray}
With arguments similar to those in {\sc Theorem \ref{th3.1}}, we
obtain
\begin{eqnarray}\label{eq4.9}
&&\|U^n_h\|^2+\sum\limits_{i=1}^n\|U^i_h-U^{i-1}_h\|^2+\frac{\Delta t \nu_0}{4} \sum\limits_{i=1}^n\|U^i_h\|^2_1\nonumber\\
&\le& C\left(\|U^0_h\|^2+\sup\limits_{1\le i\le n}\|U^i_H\|^2+\Delta t \sum\limits_{i=1}^n \|f^i\|^2 \right).
\end{eqnarray}
As $U^i_H$ satisfies inequality (\ref{eqr83.1}), we can obtain
(\ref{eq4.6}).\quad

To establish the error estimate for the solution computed by {\sc Algorithm 4.1}, we need the following lemmas.

\begin{lemma}\label{le4.3}
Let $U^n$ and $U^n_h$ be the solutions obtained by (\ref{eq2.8})-(\ref{eq2.9}) and {\sc Algorithm 4.1}, respectively. If $\Delta t$ satisfies condition \begin{eqnarray}\label{eq4.14} \Delta
t< \frac{\nu^2_0}{8\mu_1^2K^2_1T},\end{eqnarray} then for any $n\ge 1$, we
have
\begin{eqnarray}\label{eq4.15}
\frac{2}{\sqrt{\nu_0\Delta t}}\|W^n_h-W^{n-1}_h\|+\|W^n_h\|_1\lesssim H^{r}+h^{r-1}+\Delta t,
\end{eqnarray}
where $W_h^n=U^n_h-U^n$.
\end{lemma}

{\bf Proof.}
It follows from (\ref{eq2.8}) and (\ref{eq4.4}) that
\begin{eqnarray}\label{eq4.16}
\frac{1}{\Delta t}(W^n_h-W^{n-1}_h,v)+A(W^n_h,v)+\Delta t\sum\limits_{i=0}^n\omega_{ni}(\tilde B(U^i_H;U^i_h,v)-\tilde B(U^i;U^i,v))=0.~~~~~~
\end{eqnarray}
Now let us bound $|\tilde B(U^i_H;U^i_h,v)-\tilde B(U^i;U^i,v)|$. Firstly, we split $\tilde B(U^i_H;U^i_h,v)-\tilde B(U^i;U^i,v)$ as follows:
\begin{eqnarray}\label{eq4.17}&&\tilde B(U^i_H;U^i_h,v)-\tilde B(U^i;U^i,v)\nonumber\\
&=&(\alpha(U^i_H)\nabla(U^i_h-U^i),\nabla v)+((\alpha(U^i_H)-\alpha(U^i))\nabla U^i,\nabla v)+(\beta(U^i_H)-\beta(U^i),\nabla v)\nonumber\\
&&+(\gamma(U^i_H)\cdot\nabla(U^i_h-U^i),v)+((\gamma(U^i_H)-\gamma(U^i))\cdot\nabla U^i,v)+(g(U^i_H)-g(U^i),v).\nonumber\\
\end{eqnarray}
It follows that
\begin{eqnarray}\label{eq4.18}
&&|(\alpha(U^i_H)\nabla(U^i_h-U^i),\nabla v)|+|(\beta(U^i_H)-\beta(U^i),\nabla v)|\nonumber\\
&&+|(\gamma(U^i_H)\cdot\nabla(U^i_h-U^i),v)|+|(g(U^i_H)-g(U^i), v)|\nonumber\\
&\le&\mu_1\|W^i_h\|_1\|v\|_1+C_L\|U^i_H-U^i\|\|v\|_1\nonumber\\
&\le&\mu_1\|W^i_h\|_1\|v\|_1+C_L(\|u(t_i)-U^i_H\|+\|u(t_i)-U^i\|)\|v\|_1\nonumber\\
&\le&\mu_1\|W^i_h\|_1\|v\|_1+C_L(H^{r}+h^{r}+\Delta
t)\|v\|_1.\end{eqnarray} On the other hand, due to the assumption on
$\alpha$ and $\gamma$, which implies that $\alpha$ and $\gamma$ are bounded and satisfy Lipschitz condition, we have
\begin{eqnarray}\label{eq4.19}
&&\left|(\alpha(U^i_H)-\alpha(U^i))\nabla U^i,\nabla
v)\right|\nonumber\\
&\le &\left|(\alpha(U^i_H)-\alpha(U^i))\nabla (
U^i-u(t_i)),\nabla v)\right|
+\left|(\alpha(U^i_H)-\alpha(U^i))\nabla u(t_i),\nabla v)\right|\nonumber\\
&\le&  C\|\nabla (U^i-u(t_i))\|\|\nabla v\|+C_L\|U^i_H-U^i\|\|\nabla u(t_i)\|_\infty\|\nabla v\|\nonumber\\
&\le& C(u)(H^{r}+h^{r-1}+\Delta t)\|\nabla v\|,
\end{eqnarray}
and
\begin{eqnarray}\label{eq4.20}
&&\left|(\gamma(U^i_H)-\gamma(U^i))\cdot\nabla
U^i,v)\right|\nonumber\\
 &\le
&\left|(\gamma(U^i_H)-\gamma(U^i))\cdot\nabla ( U^i-u(t_i)), v)\right|
+\left|(\gamma(U^i_H)-\gamma(U^i))\cdot\nabla u(t_i),v)\right|\nonumber\\&\le&  C\|\nabla (U^i-u(t_i))\|\|v\|+C_L\|U^i_H-U^i\|\|\nabla u(t_i)\|_\infty\|v\|\nonumber\\
&\le& C(u)(H^{r}+h^{r-1}+\Delta t)\|v\|.
\end{eqnarray}
Take $v=2(W^n_h-W^{n-1}_h)$ in (\ref{eq4.16}), and combine
(\ref{eq4.16}), (\ref{eq4.17}), (\ref{eq4.18}), (\ref{eq4.19}) and
(\ref{eq4.20}) to get
\begin{eqnarray}\label{eqth421}
&&\frac{2}{\Delta t}\|W^n_h-W^{n-1}_h\|^2+\nu_0\|W^{n}_h\|^2_1+\nu_0\|W^n_h-W^{n-1}_h\|^2_1\nonumber\\
&\le& 2\Delta t
\sum\limits_{i=1}^{n}|\omega_{ni}|\left(\mu_1\|W^i_h\|_1+
C(H^{r}+h^{r-1}+\Delta t)\right)\|W^n_h-W^{n-1}_h\|_1+\nu_1\|W^{n-1}_h\|^2_1.\nonumber\\
\end{eqnarray}
The first term on the right-hand side of the above inequality can be bounded as
\begin{eqnarray}\label{eqth422}
&&2\Delta t
\sum\limits_{i=1}^{n}|\omega_{ni}|\left(\mu_1\|W^i_h\|_1+
C(H^{r}+h^{r-1}+\Delta t)\right)\|W^n_h-W^{n-1}_h\|_1\nonumber\\
&\le&\frac{4}{\nu_0}\Delta tK_1^2\mu_1^2t_n
\sum\limits_{i=1}^{n}\|W^i_h\|^2_1+
\frac{4}{\nu_0}C^2K_1^2t_n^2(H^{r}+h^{r-1}+\Delta t)^2\nonumber\\
&&+\frac{\nu_0}{2}\|W^n_h-W^{n-1}_h\|^2_1,
\end{eqnarray}
where we have used
$$2\Delta t
\sum\limits_{i=1}^{n}a_ib\le \frac{4t_n}{\nu_0}\Delta t\sum\limits_{i=1}^{n}a^2_i+\frac{\nu_0}{4t_n}\Delta t\sum\limits_{i=1}^{n}b^2=\frac{4t_n}{\nu_0}\Delta t\sum\limits_{i=1}^{n}a^2_i+\frac{\nu_0}{4}b^2.$$
Substituting (\ref{eqth422}) into (\ref{eqth421}), we
get
\begin{eqnarray}\label{eqr4.21}
&&\frac{2}{\Delta t}\|W^n_h-W^{n-1}_h\|^2+\frac{\nu_0}{2}\|W^{n}_h\|^2_1\nonumber\\
&\le&\nu_1\||W^{n-1}_h\|^2_1+\frac{4\mu^2_1K^2_1t_n}{\nu_0}\Delta t\sum\limits_{i=1}^{n}\|W^i_h\|_1+C (H^{r}+h^{r-1}+\Delta t)^2.
\end{eqnarray}
In view of (\ref{eq4.14}), application of discrete Gronwall lemma 2.1 to the above inequality yields
\begin{eqnarray}\label{eq4.22}
\frac{4}{\nu_0\Delta t}\|W^n_h-W^{n-1}_h\|^2+\|W^{n}_h\|^2_1\le C(H^{r}+h^{r-1}+\Delta t)^2 .
\end{eqnarray}
Then we arrive at (\ref{eq4.15}).\quad

Combining {\sc Theorem \ref{th3.3}} and {\sc Lemma \ref{le4.3}} immediately
yields the following theorem.
\begin{theorem}[Error estimate for two-grid FEM Algorithm 4.1]\label{th4.4}
Let $u$ be the solution of (\ref{eq2.3})-(\ref{eq2.4}) and $U^n_h$ be
the solution of Algorithm 4.1. Then, for sufficiently small $\Delta t$, we have, for all $n\ge 1$,
\begin{eqnarray}\label{eq4.23}
\|U^n_h-u(t_n)\|_1&\lesssim& H^{r}+h^{r-1}+\Delta t.
\end{eqnarray}
\end{theorem}

{\bf Proof.}
Using the triangular inequality $\|U^n_h-u(t_n)\|_1\le
\|U^n-u(t_n)\|_1+\|U^n_h-U^n\|_1$, the second inequality in
(\ref{eq3.16}), and (\ref{eq4.15}), we can obtain (\ref{eq4.23}).

From (\ref{eq4.23}), it is easy to find that
when the mesh sizes satisfy $H=O(h^{\frac{r-1}{r}})$ the two-grid
{\sc Algorithm 4.1} achieves the same approximation for PIDEs with
nonlinear memory as the classic finite element method does.

Next we will present an algorithm that reduces a nonlinear problem
to a symmetric positive definite (SPD) linear problem and a
nonlinear system of smaller size.

{\sc Algorithm 4.2.}

Step one (nonlinear problem on coarse grid $\mathcal T_H$): Given
$U^{n-1}_H$, find $U^{n}_H \in S_H$ such that
\begin{eqnarray}
  &&\frac{1}{\Delta t}\left(U^{n}_H-U^{n-1}_H,v\right)+A(U^{n}_H,v)+\Delta t\sum\limits_{i=1}^n\omega_{ni} \tilde B(U^i_H;U^i_H,v) =(f^n,v),\nonumber\\
  &&~~~~~~~~~~~~v\in S_H,~ n=1,2,\cdots,\\
&&U^0_H=u^0_H.\end{eqnarray}

Step two (SPD linear problem on fine grid $\mathcal T_h$): Given
$U^{n}_H$, find $U^{n}_h \in S_h$ such that
\begin{eqnarray}
  &&\frac{1}{\Delta t}\left(U^{n}_h-U^{n-1}_h,v\right)+A(U^{n}_h,v)+\Delta t\sum\limits_{i=1}^{n-1} \omega_{ni}\tilde B(U^i_h;U^i_h,v)
  +\Delta t\omega_{nn} \tilde B(U^n_H;U^n_H,v)\nonumber\\
   &&~~~~=(f^n,v),~~v\in S_h,\\
&&U^0_h=u^0_h,\qquad n=1,2,\cdots.\end{eqnarray}

Obviously, this algorithm can also be applied to the nonsymmetric
linear problem.
\begin{theorem}[Stability of two-grid FEM  Algorithm 4.2]
Let $U^n_h$ be the solution obtained by {\sc Algorithm 4.2}. If
$\Delta t$ satisfies (\ref{eqr3.1}), then we have
\begin{eqnarray}\label{eq4.27}
&&\|U^n_h\|+\left(\sum\limits_{i=1}^n\|U^i_h-U^{i-1}_h\|^2\right)^{1/2}+\frac{\sqrt{\nu_0}}{2}
\left(\sum\limits_{i=1}^n\Delta t\|U^i_h\|^2_1\right)^{1/2}\nonumber\\
&\le& C\left(\|U^0_h\|^2+\Delta t\|U^0_H\|^2_1+\Delta t \sum\limits_{i=1}^n \|f^i\|^2 \right)^{1/2}.
\end{eqnarray}
for any $n\ge 1$.
\end{theorem}

{\bf Proof.} Similar to (\ref{eq3.3}), using (\ref{eq4.1}), we have
\begin{eqnarray}
&&\|U^n_h\|^2-\|U^{n-1}_h\|^2+\|U^n_h-U^{n-1}_h\|^2+2\Delta t \nu_0\|U^n_h\|^2_1\nonumber\\
&\le& (\Delta t)^2
\sum\limits_{i=1}^{n-1}|\omega_{ni}|\left(\frac{\mu_0}{\epsilon}\|U^i_h\|^2_1+\mu_0\epsilon\|U^n_h\|^2_1
\right)+(\Delta t)^2|\omega_{nn}|\left(\frac{\mu_0}{\epsilon}\|U^n_H\|^2_1+\mu_0\epsilon\|U^n_h\|^2_1\right)\nonumber\\
&&+\Delta t \left(\|U^n_h\|^2+\|f^n\|^2\right).
\end{eqnarray}
After choosing $\epsilon=\frac{\nu_0}{\mu_0K_1t_n}$, the above inequality becomes
\begin{eqnarray}
&&\|U^n_h\|^2+\|U^n_h-U^{n-1}_h\|^2+\Delta t \nu_0\|U^n_h\|^2_1\nonumber\\
&\le& \|U^{n-1}_h\|^2+(\Delta t)^2
\left(\sum\limits_{i=1}^{n-1}\frac{\mu^2_0K^2_1t_n}{\nu_0}\|U^i_h\|^2_1+\frac{\mu^2_0K^2_1t_n}{\nu_0}\|U^n_H\|^2_1\right)+\Delta
t \|f^n\|^2+\Delta t \|U^n_h\|^2.\nonumber\\
\end{eqnarray}
With arguments similar to those in {\sc Theorem \ref{th3.1}}, we
obtain
\begin{eqnarray*}
&&\|U^n_h\|^2+\sum\limits_{i=1}^n\|U^i_h-U^{i-1}_h\|^2+\frac{\Delta t \nu_0}{4} \sum\limits_{i=1}^n\|U^i_h\|^2_1\nonumber\\
&\le& C\left(\|U^0_h\|^2+(\Delta t)^2\sum\limits_{i=1}^n\|U^i_H\|^2_1+\Delta t \sum\limits_{i=1}^n \|f^i\|^2 \right),\nonumber\\
\end{eqnarray*}
in view of (\ref{eqr83.1}), therefore,
\begin{eqnarray}
&&\|U^n_h\|^2+\sum\limits_{i=1}^n\|U^i_h-U^{i-1}_h\|^2+\frac{\Delta t \nu_0}{4} \sum\limits_{i=1}^n\|U^i_h\|^2_1\nonumber\\
&\le& C\left(\|U^0_h\|^2+\Delta t\|U^0_H\|^2_1+\Delta t \sum\limits_{i=1}^n \|f^i\|^2 \right),
\end{eqnarray}
which implies (\ref{eq4.27}). This completes the proof.\quad

\begin{theorem}[Error estimate for two-grid FEM Algorithm 4.2]
Let $U^n_h$ be the solution obtained by {\sc Algorithm 4.2}. Then for sufficient small $\Delta t$, we have
\begin{eqnarray}\label{eq4.31}
\|U^n_h-u(t_n)\|\lesssim \sqrt{\Delta t}H^{r-1}+h^r+\Delta t,
\end{eqnarray}
for any $n\ge 1$.
\end{theorem}

{\bf Proof.}
As in {\sc Theorem 4.4}, $W_h^n=U^n_h-U^n$ satisfies the following
error equation:
\begin{eqnarray}\label{eq4.30}
&&\frac{1}{\Delta t}(W^n_h-W^{n-1}_h,v)+A(W^n_h,v)+\Delta t\sum\limits_{i=0}^{n-1}\omega_{ni}(\tilde B(U^i_h;U^i_h,v)-\tilde B(U^i;U^i,v))\nonumber\\
&&+\Delta t\omega_{nn}(\tilde B(U^n_H;U^n_H,v)-\tilde B(U^n;U^n,v))=0.
\end{eqnarray}
In view of the assumption on the coefficients of $B$, there exists a constants $\mu_B$ such that
$$|\tilde B(u;u,v)-\tilde B(w;w,v)|\le \mu_B\|u-w\|_1\|v\|_1.$$ Then we have
\begin{eqnarray}&&|\Delta t\omega_{nn}(\tilde B(U^n_H;U^n_H,v)-\tilde B(U^n;U^n,v))|\nonumber\\
&\le&\mu_B\Delta t |\omega_{nn}|\|U^n_H-U^n\|_1\|v\|_1\nonumber\\
&\le&\mu_BK_1\Delta t (\|U^n_H-u(t_n)\|_1+\|u(t_n)-U^n\|_1)\|v\|_1\nonumber\\
&\le&\mu_BK_1\Delta t (H^{r-1}+h^{r-1}+\Delta t)\|v\|_1.
\end{eqnarray}
Take $v=2\Delta t W^n_h$ in (\ref{eq4.30}) to obtain
\begin{eqnarray}
&&\|W^n_h\|^2-\|W^{n-1}_h\|^2+\|W^n_h-W^{n-1}_h\|^2+2\Delta t \nu_0\|W^n_h\|^2_1\nonumber\\
&\le& \mu_B(\Delta t)^2 \sum\limits_{i=1}^{n-1}|\omega_{ni}|(\frac{1}{\epsilon}\|W^i_h\|^2_1+\epsilon\|W^n_h\|^2_1)\nonumber\\
&&+\mu_BK_1\Delta t \left(\epsilon\Delta t\|W^n_h\|^2_1+\frac{1}{\epsilon}\Delta t (H^{r-1}+\Delta t)^2\right).
\end{eqnarray}
By choosing $\epsilon=\frac{\nu_0}{t_n\mu_BK_1}$, we get
\begin{eqnarray}
&&\|W^n_h\|^2-\|W^{n-1}_h\|^2+\|W^n_h-W^{n-1}_h\|^2+\Delta t \nu_0\|W^n_h\|^2_1\nonumber\\
&\le& \frac{t_n\mu_B^2K_1^2}{\nu_0}(\Delta t)^2 \sum\limits_{i=1}^{n-1}\|W^i_h\|^2_1+\frac{t_n\mu_B^2K_1^2}{\nu_0}(\Delta t)^2(H^{r-1}+\Delta t)^2.
\end{eqnarray}
Sum from $1$ up to $n$ to obtain
\begin{eqnarray}
&&\|W^n_h\|^2-\|W^{0}_h\|^2+\sum\limits_{i=1}^{n}\|W^i_h-W^{i-1}_h\|^2+\Delta t \nu_0\sum\limits_{i=1}^{n}\|W^i_h\|^2_1\nonumber\\
&\le&(\Delta t)^2 \sum\limits_{i=1}^{n}\frac{t_i\mu_B^2K_1^2}{\nu_0} \sum\limits_{j=1}^{i-1}\|W^j_h\|^2_1+\sum\limits_{i=1}^{n}\frac{t_i\mu_B^2K_1^2}{\nu_0}(\Delta t)^2(H^{r-1}+\Delta t)^2\nonumber\\
&\le&(\Delta t)^2 \sum\limits_{i=0}^{n-1}\frac{t_{i+1}\mu_B^2K_1^2}{\nu_0} \sum\limits_{j=1}^{i}\|W^j_h\|^2_1+\sum\limits_{i=1}^{n}\frac{t_i\mu_B^2K_1^2}{\nu_0}(\Delta t)^2(H^{r-1}+\Delta t)^2\nonumber\\
&\le&(\Delta t)^2 \sum\limits_{i=0}^{n-1}\frac{t_{i+1}\mu_B^2K_1^2}{\nu_0} \sum\limits_{j=1}^{i}\|W^j_h\|^2_1+\frac{t_n^2\mu_B^2K_1^2}{\nu_0}\Delta t(H^{r-1}+\Delta t)^2.
\end{eqnarray}
An application of discrete Gronwall Lemma 2.1 yields
\begin{eqnarray}
&&\|W^n_h\|^2+\sum\limits_{i=1}^{n}\|W^i_h-W^{i-1}_h\|^2+\Delta t \nu_0\sum\limits_{i=1}^{n}\|W^n_h\|^2_1\nonumber\\
&\le& \frac{t^2_n\mu_B^2K_1^2}{\nu_0}\Delta t(H^{r-1}+\Delta t)^2\exp\left(\frac{t^2_n\mu_B^2K_1^2}{\nu^2_0}\right).
\end{eqnarray}
Finally, (\ref{eq4.31}) follows readily from this result when a triangular
inequality is also applied.

Next we will present an algorithm that significantly reduces
computational memory and storage requirements when $B$ gathers
lower-order spatial derivatives and nonlinear terms. To state the
algorithm, we define
$$\tilde B_s(w;u,v)=(\alpha(w)\nabla u,\nabla v),$$
and
$$N(w;u,v)=(\beta(w),\nabla v)+(\gamma(w)\cdot \nabla u+g(w),v).$$
In view of the assumptions on $\alpha(u),~\beta(u),~\gamma(u)$, and
$g(u)$, we find that there exist two constants $\mu_3$ and $\mu_4$ such that
\begin{eqnarray}\label{eqr84.1}
|\tilde B_s(w;u,v)|&\le& \mu_3\|u\|_1\|v\|_1\\
|N(w;u,v)|&\le& \mu_4\|u\|\|v\|_1.\label{eqr84.2}
\end{eqnarray}

Then the algorithm can be stated as follows.

{\sc Algorithm 4.3.}

Step one (nonlinear problem on coarse grid $\mathcal T_H$): Given
$U^{n-1}_H$, find $U^{n}_H \in S_H$ such that
\begin{eqnarray}
  &&\frac{1}{\Delta t}\left(U^{n}_H-U^{n-1}_H,v\right)+A(U^{n}_H,v)+\Delta t\sum\limits_{i=1}^n\omega_{ni} \tilde B(U^i_H;U^i_H,v) =(f^n,v),\nonumber\\
  &&~~~~~~~~~~~~v\in S_H,~ n\ge 1,\\
&&U^0_H=u^0_H.\end{eqnarray}

Step two (linear problem on fine grid $\mathcal T_h$): Given
$U^{n}_H$, find $U^{n}_h \in S_h$ such that
\begin{eqnarray}
  &&\frac{1}{\Delta t}\left(U^{n}_h-U^{n-1}_h,v\right)+A(U^{n}_h,v)+\Delta t\sum\limits_{i=1}^n\omega_{ni}
  (\tilde B_s(U^i_H;U^i_h,v)+N(U^i_H;U^i_H,v))\nonumber\\
   &&~~~~~=(f^n,v),~~~~~~~~v\in S_h,\qquad n\ge 1,\\
&&U^0_h=u^0_h,\end{eqnarray}

The stability of {\sc Algorithm 4.3} can be obtained by the same
argument for {\sc Theorem 4.1}.

\begin{theorem}[Stability of two-grid FEM  Algorithm 4.3]
Let $U^n_h$ be the solution obtained by {\sc Algorithm 4.3}. Then
when \begin{eqnarray}\label{eqrr4.6}\Delta
t\le\min\left\{\frac{1}{2},\frac{3\nu^2_0}{2\mu_3^2K^2_1T}\right\},\end{eqnarray}
we have
\begin{eqnarray}\label{eqa4.6}
&&\sup\limits_{1\le i\le
n}\|U^i_h\|+\left(\sum\limits_{i=1}^n\|U^i_h-U^{i-1}_h\|^2\right)^{1/2}+\frac{\sqrt{\nu_0}}{2}
\left(\sum\limits_{i=1}^n\Delta t\|U^i_h\|^2_1\right)^{1/2}\nonumber\\
&\le& C\left(\|U^0_h\|^2+\|U^0_H\|^2+\Delta t \sum\limits_{i=1}^n \|f^i\|^2
\right)^{1/2}.
\end{eqnarray}
\end{theorem}

{\bf Proof.} Similar to (\ref{eq3.3}), using (\ref{eqr84.1}) and (\ref{eqr84.2}), we have
\begin{eqnarray}\label{eqa4.7}
&&\|U^n_h\|^2-\|U^{n-1}_h\|^2+\|U^n_h-U^{n-1}_h\|^2+2\Delta t \nu_0\|U^n_h\|^2_1\nonumber\\
&\le& (\Delta t)^2
\sum\limits_{i=1}^{n}|\omega_{ni}|\left(\frac{\mu_3}{4\epsilon_1}\|U^i_h\|^2_1+\mu_3\epsilon_1\|U^n_h\|^2_1
+\frac{\mu_4}{4\epsilon_2}\|U^i_H\|^2+\mu_4\epsilon_2\|U^n_h\|^2_1\right)\nonumber\\
&&+\Delta t \left(\|U^n_h\|^2+\|f^n\|^2\right).
\end{eqnarray}
After choosing $\epsilon_1=\frac{\nu_0}{2\mu_3K_1t_n}$ and
$\epsilon_2=\frac{\nu_0}{2\mu_4K_1t_n}$, (\ref{eqa4.7}) becomes
\begin{eqnarray}\label{eqa4.8}
&&\|U^n_h\|^2+\|U^n_h-U^{n-1}_h\|^2+\Delta t \nu\|U^n_h\|^2_1\nonumber\\
&\le& \|U^{n-1}_h\|^2+(\Delta t)^2
\sum\limits_{i=1}^{n}\left(\frac{\mu^2_3K^2_1t_n}{2\nu_0}\|U^i_h\|^2_1+\frac{\mu^2_4K^2_1t_n}{2\nu_0}\|U^i_H\|^2\right)+\Delta
t \|f^n\|^2+\Delta t \|U^n_h\|^2.\nonumber\\
\end{eqnarray}
With arguments similar to those in {\sc Theorem \ref{th3.1}}, we
obtain
\begin{eqnarray}\label{eqa4.9}
&&\|U^n_h\|^2+\sum\limits_{i=1}^n\|U^i_h-U^{i-1}_h\|^2+\frac{\Delta t \nu_0}{4} \sum\limits_{i=1}^n\|U^i_h\|^2_1\nonumber\\
&\le& C\left(\|U^0_h\|^2+\sup\limits_{1\le i\le n}\|U^i_H\|^2+\Delta
t \sum\limits_{i=1}^n \|f^i\|^2 \right).
\end{eqnarray}
As $U^i_H$ satisfies inequality (\ref{eq3.1}), we can obtain
(\ref{eqa4.6}).\quad

To get an idea of the accuracy of {\sc Algorithm 4.3}, we have
the following theorem.
\begin{theorem}[Error estimate for two-grid FEM  Algorithm 4.3]\label{convergence}
Let $U^n_h$ be the solutions obtained by {\sc Algorithm 4.3}. Then for sufficient small $\Delta t$, we have
\begin{eqnarray}
\|U^n_h-u(t_n)\|\lesssim h^r+\Delta t+\sqrt{\Delta t}H^{r},\quad
\|U^n_h-u(t_n)\|_1\lesssim H^{r}+h^{r-1}+\Delta t,
\end{eqnarray}
 for any $n\ge 1$.
\end{theorem}

{\bf Proof.}
Set $W_h^n=U^n_h-U^n$ to get
\begin{eqnarray}
&&\frac{1}{\Delta t}(W^n_h-W^{n-1}_h,v)+A(W^n_h,v)+\Delta t\sum\limits_{i=1}^n\omega_{ni}(\tilde B_s(U^i_H;U^i_h,v)-\tilde B_s(U^i;U^i,v))\nonumber\\
&&+\Delta t\sum\limits_{i=1}^n\omega_{ni}(N(U^i_H;U^i_H,v)-N(U^i;U^i,v))=0.
\end{eqnarray}
Similar to the proof of {\sc Lemma \ref{le4.3}}, we have
\begin{eqnarray}&&N(U^i_H;U^i_H,v)-N(U^i;U^i,v)\nonumber\\
&=&(\beta(U^i_H)-\beta(U^i),\nabla v)+(\gamma(U^i_H)\cdot\nabla(U^i_H-U^i),v)\nonumber\\
&&+((\gamma(U^i_H)-\gamma(U^i))\cdot\nabla U^i,v)+(g(U^i_H)-g(U^i),\nabla v)
\end{eqnarray}
and
\begin{eqnarray}
\left|(\gamma(U^i_H)\cdot\nabla(U^i_H-U^i),v)\right|&\le& C\|U^i_H-U^i\|\|v\|_1\nonumber\\
&\le& C_L(\|u(t_i)-U^i_H\|+\|u(t_i)-U^i\|)\|v\|_1\nonumber\\
&\le& C_L(H^{r}+h^{r}+\Delta t)\|v\|_1.
\end{eqnarray}
The desired estimate can then be obtained in a way similar to proofs of {\sc Theorem \ref{th4.4}} and {\sc Lemma \ref{le4.3}}.

{\it Remark.} Observe that when $\alpha\equiv 0$, the approximation
of the integral term on the fine grid is identical to the
approximation of the integral term on the coarse grid. This means
that when we solve $U^n_h$, all $U^i_h~(i < n-1)$ do not need to be
stored on a fine grid. It also means that once the approximation of
the integral term has been computed on the coarse grid it does not
need to be computed on the fine grid. This significantly reduces
computational memory and storage requirements. This result is novel and interesting even for linear problem.

\section{Numerical experiments}

In this section, we show some experiments to confirm the effectiveness and theoretical analysis for {\sc Algorithm 4.3}. We set the domain as $[0,1]\times [0,1]$ and $T=1.00$. Noting that when $\alpha\equiv 0$ in {\sc Algorithm 4.3}, the algorithm does not need to store $U^i_h~(i < n-1)$, hence in order to confirm the efficiency and advantage of {\sc Algorithm 4.3}, we set $K(t)=e^{-t}, \alpha(u)=0, \beta(u)=(\sin u, 1-\cos u)^{T},\gamma(u)=(1-\cos u, \sin u)^{T}, g(u)=\sin u$ in \eqref{eqNV1} and we solve the following problem
\begin{equation}\label{eqNV2}
\begin{split}
&u_t-\Delta u+\int^t_{0}e^{-(t-s)}\big(-\nabla\cdot \beta(u(s))+\gamma(u(s))\cdot \nabla u+g(u(s))\big)ds\\
&~~~~~~~~~~~~~~~~~~~~~~~~~~~~~~~~~~~~~~~~~~~~~~~~~~~~~~~~~~~~~~~~~~=f(x_1,x_2; t),\\
&u(x_1,x_2; t)=0,~~~~(x_1,x_2; t)\in \partial \Omega \times (0,T],\\
&u(x_1,x_2;0)=u_0(x_1,x_2),~~~~(x_1,x_2)\in \Omega.
\end{split}
\end{equation}
We further set $u_0(x_1,x_2)=x_1(1-x_1)x_2(1-x_2)$ and
\begin{equation*}
\begin{split}
&f(x_1,x_2; t)\\
&=\big(2x_1(1-x_1)-x_1(1-x_1)x_2(1-x_2)+2x_2(1-x_2)
+(1-2x_1)x_2(1-x_2)t\big)e^{-t}\\
&-2(1-2x_1)x_2(1-x_2)e^{-t}\int_{0}^t
\cos\big(x_1(1-x_1)x_2(1-x_2)e^{-s}\big)ds\\
&+e^{-t}\int_{0}^t e^s\sin\big(x_1(1-x_1)x_2(1-x_2)e^{-s}\big)ds.
\end{split}
\end{equation*}
Then we can verify that $u(x_1,x_2;t)=x_1(1-x_1)x_2(1-x_2)e^{-t}$ is the true solution.  We use linear finite element for the space discretization. The convergence rate and effectiveness of {\sc Algorithm 4.3} in $H^1$ norm given by Theorem \ref{convergence} are confirmed in Table \ref{tab1} with $h=\frac{1}{2^l}, l=2,\cdots,9; H=\frac{1}{2}\sqrt {h}$ and $\Delta t=\frac{1}{2^l}, l=1,\cdots, 8$.

\begin{table}[h!]
\begin{center}
\begin{tabular}{|c|c|c|c|c|c|c|}
\hline
$H$&$h$&$\Delta t$&  $\|U_h^T-u(T)\|_1$ &  $h$  order& $H$ order  &$\Delta t$ order \\
\hline
$1/4$&$1/4$&$1/2$& $2.17236\times 10^{-2}$ & $--$ & $--$ & $--$\\
\hline
$1/6$&$1/8$&$1/4$ & $1.11164\times 10^{-2}$ & $0.95$ & $1.95$ & $0.95$\\
\hline
$1/8$&$1/16$&$1/8$&$5.59226\times 10^{-3}$ & $0.99$ & $1.99$ & $0.99$\\
\hline
$1/12$&$1/32$&$1/16$&$2.80089\times 10^{-3}$ & $0.99$ & $1.99$ & $0.99$\\
\hline
$1/16$&$1/64$&$1/32$&$1.40136\times 10^{-3}$ & $1.00$ & $2.00$ & $1.00$\\
\hline
$1/23$&$1/128$&$1/64$&$7.00760\times 10^{-4}$ & $0.99$ & $1.99$ & $0.99$\\
\hline
$1/32$&$1/256$&$1/128$&$3.50427\times 10^{-4}$ & $1.00$ & $2.00$ & $1.00$\\
\hline
$1/46$&$1/512$&$1/256$&$1.75207\times 10^{-4}$ & $1.00$ & $2.00$ & $1.00$\\
\hline
\end{tabular}
\end{center}
\caption{Convergence rate and accuracy of {\sc Algorithm 4.3}.}
\label{tab1}
\end{table}
Following the {\sc Algorithm 4.3}, in the numerical experiments, we do not store  $U^i_h~(i < n-1)$ and save a lot
of storege. Further, the method is much more efficient than the standard fully discrete finite element algorithm since we only need to solve a nonlinear problem with mesh-size $H=\frac{1}{2}\sqrt {h}$ and then solve the linear problem with mesh-size $h$. Using standard fully discrete finite element algorithm to solve the problem \eqref{eqNV2} by solving the nonlinear problem directly with mesh-size $h$ and the convergence rate and error in $H^1$ norm are shown in Table \ref{tab2}. Comparing Table \ref{tab1} and  Table \ref{tab2}, we can clearly see that the effectiveness and accuracy of {\sc Algorithm 4.3} are the same as standard fully discrete finite element algorithm. The error estimate in $L^2$ norm given by Theorem \ref{convergence} can also be confirmed similarly, for simplicity, we omitted listing the tables here.

\begin{table}[h!]
\begin{center}
\begin{tabular}{|c|c|c|c|c|}
\hline
$h$&$\Delta t$&  $\|U_h^T-u(T)\|_1$ &  $h$   order  &$\Delta t$ order \\
\hline
$1/4$&$1/2$& $2.17183\times 10^{-2}$ & $--$ & $--$\\
\hline
$1/8$&$1/4$ & $1.11115\times 10^{-2}$ & $0.95$& $0.95$\\
\hline
$1/16$&$1/8$&$5.58847\times 10^{-3}$ & $0.99$  & $0.99$\\
\hline
$1/32$&$1/16$&$2.79844\times 10^{-3}$ & $1.00$  & $1.00$\\
\hline
$1/64$&$1/32$&$1.39977\times 10^{-3}$ & $1.00$  & $1.00$\\
\hline
$1/128$&$1/64$&$6.99958\times 10^{-4}$ & $1.00$  & $1.00$\\
\hline
$1/256$&$1/128$&$3.49990\times 10^{-4}$ & $1.00$  & $1.00 $\\ 
\hline
$1/512$&$1/256$&$1.74996\times 10^{-4}$ & $1.00$  & $1.00$\\
\hline
\end{tabular}
\end{center}
\caption{Error and convergence rate for standard fully discrete finite element algorithm.}
\label{tab2}
\end{table}

\section{Concluding remarks}
We have presented and derived error estimates for several two-grid
finite element algorithms for PIDEs with nonlinear memory. With the backward Euler scheme, the two-grid strategy consists of
two steps: (1) discretizing the fully nonlinear problem in space on
a coarse grid with mesh-size $H$ and time step-size $\Delta t$ and
(2) discretizing the linearized problem in space on a fine grid with
mesh-size $h$ and the same time step-size as in step (1). It is
shown that these algorithms are as stable as the standard fully
discrete finite element algorithm. We also present the error
estimate at each time step. Compared with standard finite element
methods, our algorithm not only keep good accuracy but also saves a
lot of computational cost. As a byproduct of these results, we found
that one of these algorithms, {\sc Algorithm 4.3}, significantly reduces computational
memory and storage requirements if the nonlinear memory is defined
by a first-order or zero-order nonlinear differential operator. Thus, the two-grid methods studied in
this paper provide a new approach that takes advantage of some of
the nice properties hidden in a complex problem.

Numerical experiments for {\sc Algorithm 4.3} are provided to confirm the theoretical results and show that the two-grid method
has the same effectiveness and accuracy as the standard fully discrete finite element algorithm.

The analysis herein was carried out for an implicit Euler
discretization in time. However, the results could be extended to
the second-order accuracy backward differentiation formula (BDF)
scheme. Moreover, the analysis is valid for a state-dependent
forcing term $f$ that satisfies certain conditions, e.g.,
$$|\frac{\partial}{\partial u}f(x,t,u)|+|\frac{\partial^2}{\partial u^2}f(x,t,u)|\le M,~~~~ u\in \mathbb R,$$
where $M$ is a positive constant.

\section*{Acknowledgments}
The first author thanks Professor Jinchao Xu for suggesting this
problem and for many stimulating and inspiring discussions. This
paper was written at the School of Mathematical Sciences, Peking
University, where the first author spent time as a visiting
scholar.

This work was partially supported by the
National Natural Science Foundation of China [grant numers 11771060,11371074].



\end{document}